\documentclass[10pt,oneside]{article}
\usepackage[intlimits,reqno]{amsmath}
\usepackage{amsthm}
\usepackage{enumerate}
\usepackage[T1]{fontenc}
\usepackage[cp1250]{inputenc}
\usepackage{amssymb}
\newtheorem{thm}{Theorem}[section]

\newtheorem{lem}[thm]{Lemma}
\newtheorem{prop}[thm]{Proposition}
\newtheorem{rem}[thm]{Remark}
\theoremstyle{definition}

\newtheorem{I}[thm]{{\bf}}

\numberwithin{equation}{section}
\newcommand{\norm}[1]{\left\Vert#1\right\Vert}
\newcommand{\abs}[1]{\left\vert#1\right\vert}

\newcommand{\No}{\mathbb{N}_+\cup\{0\}}
\newcommand{\Z}{\mathbb Z}
\newcommand{\C}{\mathbb C}
\newcommand{\N}{\mathbb N}

\newcommand{\A}{\mathcal{A}}
\newcommand{\B}{\mathcal{B}}

\renewcommand{\L}{\mathcal{L}}

\newcommand{\Nc}{\mathcal{N}}

\renewcommand{\H}{\mathcal{H}}
\newcommand{\prf}{\paragraph{\bf Proof:}}
\renewcommand{\d}{\delta}
\renewcommand{\leq}{\leqslant}
\renewcommand{\geq}{\geqslant}
\renewcommand{\phi}{\varphi}

\newcommand{\be}{\begin{equation}}
\newcommand{\ee}{\end{equation}}
\newcommand{\bse}{\begin{subequations}}
\newcommand{\ese}{\end{subequations}}

\begin{document}

{\hfill IFT UwB / 08/ 2002}\\{\hfill August 27, 2002}

\bigskip

\begin{center}
{\bf \Large Extensions of   $C^*-$algebras by partial isometries}\\
\bigskip
{\large A. Lebedev,   A. Odzijewicz\footnote{This work is supported in part by KBN grant 2 PO3 A 012 19.}  }\\ \bigskip
  Institute of Mathematics,  University  in Bialystok,
 ul. Akademicka 2, PL-15-424 Bialystok, Poland\\

Institute of Theoretical Physics,  University  in Bialystok,
 ul. Lipowa 41, PL-15-424 Bialystok, Poland
\end{center}
\bigskip
\begin{abstract}
In the present paper we study the structure of $C^*-$algebras generated by a certain
$^*-$algebra $\A$ and a partial isometry inducing an endomorphism of $\A$.
\end{abstract}



\tableofcontents

\bibliographystyle{abbrv}

%\title{Factorization method for functional equations of second order}
%\author{Anatol Odzijewicz\thanks{aodzijew@labfiz.uwb.edu.pl}, Tomasz Goliński\thanks{tomaszg@alpha.uwb.edu.pl}}
\section{Introduction}

The $C^*-$algebras associated to automorphisms (endomorphisms) play a significant role as in
different fields
of mathematics
 so also in mathematical physics (see, for example,
  \cite{David,Peder,AntLeb,AnLebBel1,AnLebBel2,OstSam})
and in particular in quantum optics (see \cite{OdHorTer,OHT2}). At present the theory of
$C^*-$algebras
associated to automorphisms is extremely well developed (recall again \cite{Peder}).
On the other hand the theory of $C^*-$algebras associated to endomorphisms has not
attained a complete shape yet though
a progress here is rather rapid.
The main achievements in this direction have been obtained under the assumption
that an endomorphism of a $C^*-$algebra is generated by an {\em isometry} (see, for example,
\cite{BKR,ALNR} and the sources cited there). But in fact
the most natural situation is that when an endomorphism is generated by a {\em partial isometry}
and namely the $C^*-$algebras associated to the endomorphisms of this type are the subject of this
paper.

The article is organized as follows.

In the second section (first after introduction) we give the principal results describing the
structure of the $C^*-$algebra $\B = \B(\A, U)$ generated by a $^*-$algebra $\A \subset L(H)$
and a partial isometry $U\in L(H)$ such that the mapping
$$
\A \ni a \mapsto UaU^*
$$
is an endomorphism of $\A$. This description is carried out under
 the assumption that $\A$ and $U$ satisfy the following three conditions
 $$
Ua = UaU^*U, \ \ \ a\in \A ;
 $$
 $$
UaU^* \in \A , \ \ \ a\in \A
 $$
 and
 $$
U^*aU \in \A , \ \ \ a\in \A .
 $$
We call the algebras $\A$ possessing these properties  the {\em coefficient algebras}
 (for $\B$).

 It is shown that any element of $\B$ can be presented as a Fourier like series with coefficients
 from $\A$ (Theorems \ref{uniq}, \ref{uniqueNk}). In addition the isomorphism theorem
  (Theorem \ref{iso}) uncovering the uniqueness of the 'internal' structure of $\B$ in the
  presence of a certain property $(*)$ (see \ref{*'}) is proved.

  The third section gives a description of the construction of  coefficient
   algebras starting from certain initial algebras that are not coefficient algebras.
   Further in the  fourth section we give the corresponding
   construction of {\em commutative} coefficient algebras.

   Finally in the last section we present two examples of $C^*-$algebras of the type
   considered. The first one is related to the polar decomposition of an operator and the second
   arise from the problems of quantum optics.

   In connection with the subject of the paper we should also mention the article
   \cite{LebOdz1}
   where the subcase of the algebra considered here related to the polar decomposition of an
   operator is investigated (see also Example 1 in Section 5).

\section{Extensions by partial isometries. Isomorphism theorem}

Throughout this section we fix a $^*$-algebra $\A\subset L(H)$ containing the identity $1$ of $L(H)$
and proceed to the description of the $C^*-$extensions of $\A$ associated with the mappings
\begin{equation}
\d(x)=UxU^*,\qquad \d_*(x)=U^*xU,\qquad x\in L(H)
\end{equation}
where $U\in L(H), U\neq 0$. It is clear that
$\d$ and $\d_*$ are linear and continuous ($\norm\d=\norm{\d_*}=\norm U^2$) maps of $L(H)$ and
$\d(x^*)=\d(x)^*$, \ \ $\d_*(x^*)=\d_*(x)^*$. When using the powers $\d^k$ and $\d_*^k,\ \
k=0,1,2,\ldots $ we
 assume for
convenience that $\d^0 (x) = \d_*^0 (x) = x $.

Observe that if $\d : \A \to L(H)$ is a morphism then we have
 $$
UU^* = \d (1) = \d (1^2) = \d^2 (1) = (UU^*)^2
  $$
and therefore $UU^*$ is a partial isometry. Because of this partial isometries
play a principal role in the subject of the paper.

{\bf Remark}. If $U$ is an {\em isometry} (that is $U^*U =1$) then for {\em  any} algebra
$\A$\ \
 $\d : \A \to
L(H)$ is a morphism  since for any  pair of elements $a,b \in \A$ we have
 $$
 \d (ab) = UabU^*  = UaU^*UbU^* = \d (a)\d (b)
 $$

We recall that a
linear bounded operator $U$ in a hilbert space $H$ is called a
partial isometry if there exists a closed subspace $H_1 \subset H$
such that $$ \Vert U\xi \Vert = \Vert \xi  \Vert , \ \ \ \xi \in
H_1 $$ and $$ U\xi = 0, \ \ \ \ \xi \in H\ominus H_1. $$ The space
$H_1$ is called the {\em initial} space of $U$ and $U(H_1)$ is
called the {\em final} space of $U$.

\begin{rem}
\label{2.1}
\rm
 Hereafter we list the  well known equivalent characteristic
properties of a partial isometry:\\

\noindent 1) $U$ is a partial isometry,\\
2) $U^*$ is a partial isometry,\\
3) $U^*U$ is a projection (onto the initial space of $U$),\\
4) $UU^*$ is a projection (onto the final space of $U$),\\
5) $UU^*U =U$ and \ $U^*UU^* =U^*$.\\
(see for example \cite{Halm}, problem 98).
\end{rem}

We shall study the $C^*$-algebra $\B:=\B(\A,U)$ generated by $\A$ and $U$ assuming additionally
that $\A$ is the \textbf{coefficient algebra} of $\B$,
by this we mean  that  $\A$ possesses the following three
properties
\begin{equation}
\label{1}
 Ua=\d(a)U, \ \ \ \ a\in \A ,
\end{equation}
 \begin{equation}
\label{one}
 \d:\A\to\A ,
 \end{equation}
 \begin{equation}
\label{two}
 \d_*:\A\to\A .
 \end{equation}
 The algebras possessing these properties really play the role of the
'coefficients' in $\B(\A,U)$ which is shown in Proposition \ref{B0}. In the next Section 3 we
shall discuss the construction of the algebra satisfying \eqref{1}, \eqref{one} and \eqref{two}
starting from an initial algebra that satisfies only some of these conditions
or even does not satisfy any of them.

Note that property \eqref{1} is equivalent to the property
\begin{equation}
\label{1'}
 aU^*=U^*\d(a), \ \ \ \ a\in \A
\end{equation}
 which can be verified by passage to the adjoint operators.

 The next useful observation shows that property \eqref{1} can be also written in a different way.

\begin{prop}
\label{equiv}
 Let $\A$ be a subalgebra of $L(H)$, $1\in\A$ and $U\in L(H)$ then
the following three conditions are equivalent\\

(i)\ \ $\A$ and $U$ satisfy  condition \eqref{1};\\

 (ii)\ \ U is a partial isometry and \begin{equation}
\label{comm}U^*U\in\A',\end{equation} where $\A'$ is the commutant of $\A$;\\

(iii) \ \   $U^*U\in\A'$ and $\d : \A \to \d(\A )$ is a morphism.
\end{prop}
\prf \ \ (i)\ \ $\Rightarrow$\ \ (ii). Let \eqref{1} and so \eqref{1'} be satisfied. Taking $a=1$
in  \eqref{1} one has $U= \d (1)U= UU^*U$. Thus $U$ is a partial isometry.

Multiplying \eqref{1} by $U^*$ from the left and applying \eqref{1'} we obtain
 $$
 U^*Ua =U^*\d (a)U= aU^*U.
 $$
So $U^*U \in A'$.\\

(ii)\ \ $\Rightarrow$\ \ (i). Since $U$ is a partial isometry and  $U^*U\in A'$ we have for any
$a\in \A$
 $$ Ua = UU^*Ua =UaU^*U =\d (a)U . $$
 Thus \eqref{1} is true.\\

(ii)\ \ $\Rightarrow $\ \ (iii). \ \ It is enough to
 show that
 $$\d (ab)= \d (a)\d (b), \ \ a,b\in \A .$$
  But this follows from the fact that $U$ is a partial isometry $U^*U \in \A'$ and the
relations
 $$
  \d(ab) = UabU^* =UU^*U abU^* = UaU^*UbU^* =\d (a)\d (b).  $$

(iii)\ \ $\Rightarrow $\ \ (ii). \ \ Since $\d$ is a morphism and $1\in \A$ it follows that $$
UU^*= \d (1)  = \d(1^2) = \d(1)\d(1) = (UU^*)^2.
 $$
Thus $U$ is a partial isometry.  \qed \\

From this proposition it follows that in the above  definition of the
  coefficient algebra one can replace   \eqref{1}, for example,
  by \eqref{comm} and assumption that $U$
  is a partial isometry.

Now we present the first important result on the structure of the algebra
$\B$.
\begin{prop}
\label{B0} Let $\A$ and $U$ satisfy  conditions \eqref{1},  \eqref{one} and
\eqref{two}.
 Then the
vector space $\B_0$ consisting of finite sums
\begin{equation}
\label{suma}
x= U^{*N}a_{\bar N}+ \ldots +
U^*a_{\bar 1}+ a_0 + a_1U +\ldots + a_NU^N ,
\end{equation}
 where $a_k, a_{\bar l}\in\A$ and $N\in\N\cup\{0\}$, is a uniformly
dense $^*-$subalgebra of the $C^*$-algebra $\B$.
\end{prop}
\prf Clearly if $b\in \B_0$ then $b^* \in \B_0$. Let us verify that $\B_0$ is an algebra.

Routine computation using \eqref{1},   \eqref{one} and \eqref{two} shows that for the product
$$aU^kU^{*l}b, \ \ a,b\in \A $$ we have the following possibilities:\\

1) \ \ $k\leq l$ : $aU^kU^{*l}b = U^{*l-k}\d^{l-k} (a)\d^l (1)b\in U^{*l-k}\A$;

2) \ \ $k\geq l$ : $aU^kU^{*l}b = a\d^k (1)\d^{k-l}(b)U^{k-l} \in \A U^{k-l}$.\\
 And for the
product $$U^kaU^{*l}, \ \ a\in \A $$ we have the following possibilities:\\

1)\ \ $k\leq l$ : $U^kaU^{*l}= U^{*l-k}\d^l (a)\in U^{*l-k}\A$;

2)\ \ $k\geq l$ : $U^kaU^{*l}= \d^k (a)U^{k-l}\in \A U^{k-l}$.

These relations along with the   properties  \eqref{1},   \eqref{one} and \eqref{two} in turn
imply that $\B_0$ is an algebra.

Now, let $\B_f$ be the algebra of all finite algebraic combinations of $U^k$, $U^{*l}$, $k,l\in\N$
and elements of the algebra $\A$. Clearly $$\overline{\B_f} =\B .$$ But the above relations also
imply the equality $\B_f =\B_0$ which finishes the proof.
 \qed
 \begin{rem}
\rm If one assumes that additionally to the assumptions of Proposition\ref{B0} we have
 \begin{equation}
\label{01'}
 U^*a=\d_*(a)U^*\;\Leftrightarrow\;aU=U\d_*(a)
  \end{equation}
 for $a\in\A$, that is if $\d_*$  also satisfies
 condition \eqref{1} then the statement proved above will be valid for the subalgebra $\B_0$
generated by the elements
\begin{equation}
 x=a_{\bar N}U^{*N}+\ldots + a_{\bar 1}U^* +a_0+ a_1U+\ldots + a_NU^N .
\end{equation}
This modification is possible since by \eqref{01'} one can put the coefficients $a_{\bar k}$ in
front of the operators $U^{*k}$.
\end{rem}
 Let us observe that using $U^k=U^kU^{*k}U^k$ and $U^{*k}=U^{*k}U^kU^{*k}$
(see Proposition \ref{Uk})
 one can always choose
coefficients $a_k$ and $a_{\bar l}$ of \eqref{suma} in such a way that
\begin{equation}
\label{wybor}
a_kU^kU^{*k}=a_k\qquad\textrm{and }\qquad U^kU^{*k}a_{\bar k}=a_{\bar k}.
\end{equation}
%\begin{equation} x=(a_NU^NU^{*N})U^N+\ldots+a_1U+a_0+U^*a_{\bar 1}+\ldots+U^{*N}a_{\bar N} \end{equation}|
However the assumption that $a_k,a_{\bar k}\in\A$, $k=1,\ldots,N$ satisfy \eqref{wybor} does not
guarantee yet their uniqueness in the expansion \eqref{suma}. The uniqueness of the coefficients
will be achieved further in Theorem \ref{uniq} and Theorem \ref{uniqueNk}
in  the presence of the next property of $\B_0$.
\begin{I}
\label{*'}
{\rm \,
We shall say that the algebra $\B_0$ possesses the\ \ {\em property} (*)\ \ if for any $x\in\B_0$
given by \eqref{suma} the inequality \begin{equation} \label{gwiazdka}\norm{a_0}\leq \norm{x} \qquad\qquad (*)
\end{equation} holds.}
\end{I}
\begin{prop}
\label{gwiazdka2p}
 Let $\B_0$ be the algebra considered in Proposition \ref{B0}.
 If $\B_0$ possesses  property (*) and coefficients of $x\in\B_0$ do satisfy
\eqref{wybor} then
\begin{equation}
\label{gwiazdka2}
\norm{a_k}\leq \norm x
\end{equation}
$$
\norm{a_{\bar k}}\leq \norm{x}
 $$
 for $k,l\in\{0,1,\ldots,N\}$.
\end{prop}
\prf By means of routine computation one can show  that $0$-degree components of elements $U^kx$
%and $U^{*k}x$
are
% \begin{equation}
% \label{xU*}
% (xU^{*k})_0=a_kU^kU^{*k}
 %\end{equation}
%  and
 \begin{equation}
 \label{Ux}
   (U^kx)_0= U^kU^{*k}a_{\bar k}.  \end{equation} \\
 Since $ \Vert U^k \Vert \leq 1$ it follows that
\begin{equation}
\label{>}
 \Vert x \Vert \geq \Vert U^k x \Vert
 \end{equation}
Now \eqref{Ux}, \eqref{>} and property (*) imply
\begin{equation}
\label{*'}
 \Vert x \Vert \geq \Vert U^k x \Vert \geq \Vert (U^kx)_0  \Vert =
   \Vert U^kU^{*k}a_{\bar k}\Vert = \Vert a_{\bar k} \Vert .
 \end{equation}
 Since $\Vert x\Vert = \Vert x^* \Vert$ inequality \eqref{*'} being applied to $x^*$ gives
 $$
\Vert x\Vert = \Vert x^* \Vert \geq \Vert a_k^* \Vert =\Vert a_k \Vert
 $$
 \qed

Now Propositions \ref{B0} and  \ref{gwiazdka2p} give the following uniqueness result.
\begin{thm}
\label{uniq} Let $\A$ and $U$ satisfy  conditions \eqref{1} (or \eqref{comm}),  \eqref{one} and
\eqref{two}.
 Then property (*) implies the
uniqueness of the decomposition \eqref{suma} of any element $x\in\B_0$ with the coefficients
satisfying \eqref{wybor}.
\end{thm}
 Hence in the presence of the property (*) one can define the linear and continuous
maps $\Nc_k : \B_0 \to \A$ and $\Nc_{\bar k} : \B_0 \to \A , \ \ k\in\No$
\begin{equation}
\label{N-k} \Nc_k(x)=a_k\in\A U^kU^{*k}\subset\A \end{equation}
 $$
 \Nc_{\bar k}(x)=a_{\bar k}\in U^kU^{*k}\A\subset\A.
 $$
By continuity these mappings can be expanded onto the whole of $\B$ thus defining the {\bf
'coefficients'} of an arbitrary element $x\in \B$. We shall show further in Theorem \ref{uniqueNk}
  that these
coefficients determine $x$ in a unique way.\\

Theorem \ref{3a.N}  presented  below shows that once the algebra $\B_0$ possesses the property (*)
the norm of an element $x\in \B_0$ can be calculated only in terms of the elements of $\A$
(0-degree coefficients of the powers of $xx^*$). Among the key moments in the proof of this
theorem are the norm estimates of sums of elements in $C^*-$algebras listed in the next Lemma
\ref{3a.sum}.  The estimates presented in this lemma are useful in their own right and probably
 known (in particular the components of the statement of the  lemma are given in
\cite{Ant}, Lemma 7.3 and \cite{AntLeb}, Lemma 22.3)). The proof of the lemma can be obtained as a
simple modification of the reasoning given in the proof of \cite{AntLeb}, Lemma 22.3.\\
\begin{lem}
\label{3a.sum} For any $C^*-$ algebra $B$ and any elements $d_1, ..., d_m \in B$ we have
\begin{equation}
\label{be3.81} \left\Vert  \sum_{i=1}^m d_i  \right\Vert^2 \le m \left\Vert  \sum_{i=1}^m d_i
d_i^* \right\Vert
\end{equation}
and
\begin{equation}
\label{be3.81a} \left\Vert  \sum_{i=1}^m d_i  \right\Vert^2 \le m \left\Vert  \sum_{i=1}^m d_i^*
d_i \right\Vert
\end{equation}
On the other hand
\begin{equation}
\label{be3.82} \left\Vert  \sum_{i=1}^m |d_i|  \right\Vert^2 \ge\frac{1}{m} \left\Vert
\sum_{i=1}^m d_i^* d_i \right\Vert
\end{equation}
and
\begin{equation}
\label{be3.82a} \left\Vert  \sum_{i=1}^m \sqrt{d_i d_i^*}  \right\Vert^2 \ge\frac{1}{m} \left\Vert
\sum_{i=1}^m d_i d_i^* \right\Vert
\end{equation}
\end{lem}
\begin{thm}
\label{3a.N} Let   $\A$ and $U$ satisfy the assumptions of Theorem \ref{uniq}. If $\B_0$ possesses
the property (*) then for any element $x$ of the form  \eqref{suma} we have
\begin{equation}
\label{be3.131} \Vert x \Vert = \lim_{k\to\infty} \sqrt[4k]{ \left\Vert \Nc_0 \left[
(xx^*)^{2k}\right]\right\Vert }
\end{equation}
where $\Nc_0$ is the mapping defined by \eqref{N-k}
\end{thm}
\prf Applying (\ref{be3.81}) to the operator $$ x= U^{*N}a_{\bar N}+ \ldots + U^*a_{\bar 1}+ a_0 +
a_1U +\ldots + a_NU^N = d_{\bar N} + ... + d_0 + ... + d_N $$
 we obtain
 $$
 \Vert x \Vert^2 \leq
(2N+1) \left\Vert \sum_{i=0}^N d_i d_i^* +  \sum_{i=1}^N d_{\bar i} d_{\bar i}^* \right\Vert =
(2N+1) \Vert \Nc_0 (xx^*)\Vert $$
 where
 $$
 d_id_i^* =
a_iU^iU^{i*}a_i^*, \ \ i=0, \ldots , N
 $$
  $$
   d_{\bar i} d_{\bar i}^* = U^{*i}a_{\bar i}a_{\bar i}^* U^i , \ \ i=1,\ldots , N
 $$
  and therefore  $d_id_i^*, d_{\bar i} d_{\bar i}^* \in  \A $.\\
 On the other hand as
$\B_0$ possesses the property $(*)$ we have
$$
 \Vert x \Vert^2 = \Vert xx^* \Vert \ge \Vert \Nc_0
(xx^*)\Vert
$$
 thus
\begin{equation}
\label{be3.101} \Vert \Nc_0 (xx^*)\Vert \le  \Vert xx^* \Vert = \Vert x \Vert^2  \le (2N+1) \Vert
\Nc_0 (xx^*)\Vert
\end{equation}
Applying (\ref{be3.101}) to $(xx^* )^k$ and having in mind that $(xx^* )^k = (xx^* )^{k*} $ and
$\Vert (xx^* )^{2k}\Vert = \Vert x \Vert^{4k} $ one has
$$
\Vert \Nc_0 \left[
(xx^*)^{2k}\right]\Vert \le
 \Vert (xx^*)^k \cdot  (xx^*)^{k*}  \Vert =
\Vert x \Vert^{4k}  \le (4kN+1) \Vert \Nc_0 \left[ (xx^*)^{2k}\right]\Vert
$$
 since being written
in the form \eqref{suma} $(xx^*)^k $ has not more than $(4kN+1)$ summands.\\ So $$
 \sqrt[4k]{ \left\Vert \Nc_0 \left[ (xx^*)^{2k}\right]\right\Vert }\le \Vert x \Vert \le
\sqrt[4k]{ 4kN+1 }\cdot \sqrt[4k]{  \left\Vert \Nc_0 \left[ (xx^*)^{2k}\right]\right\Vert }
 $$
Observing the equality $$
 \lim_{k\to\infty}\sqrt[4k]{ 4kN+1 } =1
 $$
  we conclude that
  $$
   \Vert x
\Vert = \lim_{k\to\infty} \sqrt[4k]{ \left\Vert \Nc_0 \left[ (xx^*)^{2k}\right]\right\Vert }.
 $$
  The
proof is complete. \qed
\begin{rem}
\rm
The formulae of \eqref{be3.131} type have been the subject of study in connection with various
algebraic and invertibility problems of the theory of operators. The main ideas concerning the
spectral radius evaluation of elements of $C^*-$algebras associated to automorphisms by means of
the formulae of this type arise in the work by  Brenner \cite{Bren1,Bren2,Bren3,Bren4}
 where the corresponding
results were obtained in the case of a commutative algebra $\A$ and a unitary representation of
the group $\Z \to L(H)$ inducing automorphisms of $\A$. It has been also observed  there
that these formulae can be generalized to the case of a subexponential group.  Using these
formulae  Brenner has found a proof of the isomorphism theorem (of Theorem \ref{iso} type)
and this method was generalized further in  \cite{Ant11,Ant}. The general
'noncommutative' formulae of this type for the evaluation of the spectral radius are given in
\cite{AntLeb}, Section 22. It is worth mentioning that perhaps the first spectral radius formulae
of this type appeared in the classical paper by Beurling \cite{Beurl}
\end{rem}

Now we are ready to formulate the main result of this section
 which is a generalization of the  {\bf
isomorphism theorem} for the group ($\Z,+)$ (see \cite{AntLeb}, Section 12) to the case when one replaces
automorphisms related to unitary operators by endomorphism related to a partial isometry.

\begin{thm}
\label{iso}
Let $(\A_i,U_i), \ \ i=1,2$ satisfy  conditions \eqref{1}, \eqref{one}, \eqref{two} and
$\B_{0i}$ possess  property (*).\\
 Suppose that
\begin{equation}
\label{A1-A2}
 \phi:\A_1\to \A_2 \end{equation}
 is an isomorphism and \begin{equation}
\phi\circ\d_1=\d_2\circ\phi.\end{equation} Then the map \begin{equation} \Phi(x):=\phi(x),\qquad x\in\A_1\end{equation} $$
\Phi(U_1):=U_2$$ gives rise to the isomorphism between the $C^*$-algebras $\B_1 =\B (\A_1,U_1)$
and $\B_2 = \B (\A_2,U_2)$.
\end{thm}
\prf  Consider an operator\ \ $x \in  \B_{01}$\ \ having the form
\begin{equation}
\label{e3.01} x  = U_1^{N*}a_{\bar N} + . . . +a_0 + . . . + a_N U_1^N ,
\end{equation}
where $a_{\bar k}, a_k \in\,  \A_1, \ \ k=0,..., N $.\\
 Let
 \ \ $\Phi (x) \in \B_{02}$ be the operator given by
\begin{equation}
\label{e3.02} \Phi (x)  = U_2^{N*}\phi (a_{\bar N}) + . . . +\phi ( a_0) + . . .
 + \phi (a_N)U_2^N .
\end{equation}
The assumptions of the theorem imply that $\Phi $ establishes a *-algebraic isomorphism between
$\B_{01}$ and $\B_{02}$ therefore to finish the proof  it is enough to verify the equality\
\
 $\Vert x \Vert = \Vert \Phi (x)
\Vert , \ \ \ \ \ x\in \B_{01}$.\\
 By Theorem \ref{3a.N} we have
\begin{equation}
\label{be3.013} \Vert x \Vert = \lim_{k\to\infty} \sqrt[4k]{ \left\Vert \Nc_0 \left[
(xx^*)^{2k}\right]\right\Vert }
\end{equation}
where $$
 \Nc_0 :  \B_{01} \to\,  \A_1
$$ is described in \eqref{N-k}.\\ Similarly
\begin{equation}
\label{be3.013'} \Vert \Phi (x) \Vert = \lim_{k\to\infty} \sqrt[4k]{ \left\Vert \Nc_0 \left[ (\Phi
(x) \Phi (x^*))^{2k}\right]\right\Vert }
\end{equation}
where
 $$
 \Nc_0 :  \B_{02} \to \A_2 .
$$
Observe that (\ref{e3.01}) and (\ref{e3.02}) along with the assumptions of the theorem imply
 $$
 \Nc_0 \left[ (\Phi (x) \Phi (x^*))^{2k}\right] =  \Nc_0 \left[ (\Phi (xx^*))^{2k}\right] =
\phi \left(  \Nc_0 \left[ (xx^*)^{2k}\right] \right)
 $$
  and therefore (in view of \eqref{A1-A2}))
$$
 \left\Vert \Nc_0 \left[ (\Phi (x) \Phi (x^*))^{2k}\right] \right\Vert  =
 \left\Vert  \Nc_0 \left[ (xx^*)^{2k}\right] \right\Vert
$$
 This along with (\ref{be3.013'}) and  (\ref{be3.013}) implies the equality
 $$
  \Vert x \Vert =
\Vert \Phi (b) \Vert $$ and finishes the proof.\qed
\begin{rem}
\label{cross} \rm It is worth mentioning that   property $(*)$ and the results of Theorem
\ref{iso} type play a fundamental role in the theory of crossed products of $C^*-$algebras by
discrete groups (semigroups) of automorphisms (endomorphisms). Namely this property is a {\em
characteristic} property of the crossed product and it enables one to construct its  faithful
representations. The importance of  property $(*)$ for the first time (probably) was clarified
by O'Donovan \cite{O'Donov} in connection with the description of  $C^*-$algebras generated by
weighted shifts. The most general result establishing the crucial role of this property in the
theory of crossed products of $C^*-$algebras by discrete groups of {\em automorphisms} was
obtained in \cite{Leb1} (see also \cite{AntLeb}, Chapters 2,3 for complete proofs and various
applications) for an {\em arbitrary} $C^*-$algebra and {\em amenable} discrete group. The relation
of the corresponding property to the faithful representations of crossed products by {\em
endomorphisms} generated by {\em isometries} was investigated in \cite{BKR,ALNR}. The properties
of this sort proved to be of great value not only in pure $C^*-$theory but also in various
applications such as, for example, the construction of  symbolic calculus and developing the
solvability theory of functional differential equations (see \cite{AnLebBel1, AnLebBel2}).
\end{rem}

We finish the section with the statement showing that in the presence of property $(*)$
 any element $x\in \B$ can be 'restored' by its
 coefficients $\Nc_k (x)$ and $\Nc_{\bar k} (x)$.

\begin{thm}
\label{uniqueNk}
Let $\B_0$ possess property $(*)$. If $x\in \B$ is such that $\Nc_k (x)= 0,\ \ k = 0,1,\ldots $
and $\Nc_{\bar{k}} (x) =0, \ \ k = 1,\ldots$ then $x=0$.
\end{thm}
\prf Let $S^1 = \{ \lambda \in \C : \abs \lambda =1 \}$ be the unite circle on the complex plane.
For every  $\lambda \in S^1$ we consider
 the algebra $\B_\lambda := \B (\A , \lambda U)$ (clearly $\B =\B_\lambda$). Since $\B_0$ possesses property $(*)$
 the algebra $\B_{\lambda 0}:= (\B_\lambda )_0$ possesses property $(*)$ as well and
 Theorem \ref{iso} tells us that the mapping
$$
\Phi_\lambda (a)= a, \ \ a\in \A; \ \ \ \Phi (U) = \lambda U \
$$
establishes an isomorphism between $\B$ and $\B_\lambda$.
 To shorten the notation we set for every
$x\in \B$ and  $\lambda \in S^1$
$$
x(\lambda ) = \Phi_\lambda (x).
$$
The mentioned isomorphism gives
$$
\Vert x \Vert =
\Vert x (\lambda ) \Vert , \ \ \ x\in \B
$$
and in particular for every finite sum
$$ y= U^{*N}a_{\bar N}+ \ldots + U^*a_{\bar 1}+ a_0 +
a_1U +\ldots + a_NU^N
 $$
with $a_{\bar k}, a_k $ satisfying \eqref{wybor} we have
\begin{equation}
\label{ravno}
\Vert   y \Vert = \Vert   y(\lambda) \Vert =
\Vert  \lambda^{-N} U^{*N}a_{\bar N}+ \ldots +  \lambda^{-1} U^*a_{\bar 1}+ a_0 +
 \lambda a_1U +\ldots + \lambda ^N a_NU^N \Vert
\end{equation}
To prove that $x=0$ it is enough to show that for any fixed $\xi,\eta \in H$
with $\Vert \xi \Vert = \Vert \eta \Vert = 1$ we have
\begin{equation}
\label{zero}
<x\xi , \eta >= 0
\end{equation}
where $<\, ,\, >$ is the inner product in $H$.

Let $x_n , \ \ n=1,2, \ldots $ be a sequence of elements of $\B$ tending to $x$ and each
having the form
 $$
 x_n= U^{*N}a^{(n)}_{\bar {N}}+ \ldots + U^*a^{(n)}_{\bar {1}}+ a^{(n)}_0 +
a^{(n)}_1U +\ldots + a^{(n)}_NU^N
$$
where $N = N(n)$ and $a^{(n)}_{\bar {k}}, a^{(n)}_k,\ \ k=0,1, \ldots , N(n)$ satisfy
\eqref{wybor}. Consider the elements $x_n (\lambda )$ and the sequence
$f_n (\lambda ), \ \ n=1,2, \ldots$
of continuous functions on $S^1$ defined by
\begin{equation}
\label{flambda}
f_n (\lambda ) = < x_n (\lambda )\xi, \eta > = \sum_{k=- N(n)}^{N(n)}\alpha_k^{n}\lambda^k
\end{equation}
where
$$
\alpha_k^{n} = <a_k^{n}U^k \xi, \eta>, \  \  k=0,1,\ldots ,N(n)
$$
and
$$
\alpha_{-k}^{n} = < U^{*k}a^{(n)}_{\bar {k}}\xi ,\eta >, \ \ k=1,2,\ldots N(n).
$$
Since $x_n \to_{n\to \infty} x$ it follows (in view of \eqref{ravno} ) that
\begin{equation}
\label{xn>}
\Vert x_{n_1} (\lambda )- x_{n_2}(\lambda )\Vert =
\Vert x_{n_1} - x_{n_2}\Vert \to_{n_1, n_2 \to\infty}0 .
\end{equation}
Therefore for every fixed $\lambda_0 \in S^1$ the sequence
$x_n (\lambda_0)$ tends to a certain element $x(\lambda )\in \B$.

Let $f(\lambda )$ be the function given by
$$
f(\lambda )= <x(\lambda )\xi ,\eta >.
$$
Applying \eqref{xn>} we obtain
$$
\abs{f_{n_1} (\lambda )- f_{n_2}(\lambda )} =
\abs{<[x_{n_1} (\lambda )- x_{n_2}(\lambda ) ]\xi ,\eta >}\le
\Vert x_{n_1} - x_{n_2}\Vert \to_{n_1, n_2 \to\infty}0
$$
which means that the sequence $f_n$ of (continuous) functions
tends uniformly (on $S^1$) to $f$. Thus $f$ is continuous and
therefore $f\in L^2 (S^1)$.

Let
$$
f (\lambda ) = \sum_{k=-\infty }^\infty \alpha_k \lambda ^k
$$
where the righthand part is the Fourier series of $f$. Since $f_n \to f$
(in $L^2 (S^1)$) it follows that
\begin{equation}
\label{>10}
\alpha^{(n)}_k \to \alpha_k\ \ {\rm  for\ \ every }\ \     k
\end{equation}
where $\alpha^{(n)}_k$ are those defined by \eqref{flambda}.

Now note that property $(*)$ implies
\begin{equation}
\label{>11}
\Vert a^{(n)}_k \Vert \to \Vert \Nc_k (x)   \Vert \ \ {\rm  for\ \ every}\ \ k\ge 0
\end{equation}
and
\begin{equation}
\label{>11'}
\Vert a^{(n)}_{\bar k} \Vert \to \Vert \Nc_{\bar k} (x)   \Vert \ \ {\rm  for\ \ every}\ \ k\ge 0
\end{equation}
And also observe that
$$
\abs {\alpha^{(n)}_k} \le \Vert a^{(n)}_k \Vert \ \ {\rm  for\ \ every}\ \ k\ge 0
$$
and
$$
\abs {\alpha^{(n)}_{-k}} \le \Vert a^{(n)}_{\bar k} \Vert\ \ {\rm  for\ \ every}\ \ k\ge 0
$$
which together with \eqref{>10}, \eqref{>11} and \eqref{>11'}
means that
\begin{equation}
\label{>12}
\alpha_k = 0 \ \ {\rm  for\ \ every}\ \ k \in \Z .
\end{equation}
Now \eqref{>12} and the continuity of $f$ implies
$$
f(\lambda )=0  \ \ {\rm  for\ \ every}\ \ \lambda \in S^1 .
$$
In particular
$$
f(1) =<x\xi , \eta >=0 .
$$
Thus \eqref{zero} is true and the proof is finished. \qed

\section{Coefficient algebra}

For the constructions given in the previous section it is important to have a $^*-$algebra $\A$
satisfying the conditions \eqref{1} (or \eqref{comm}),  \eqref{one} and \eqref{two} that is a
coefficient algebra. This section is devoted to the investigation of the method of constructing
  the coefficient algebra starting from a suitable initial algebra. The principal result here is
  Theorem \ref{coef-ext}.

Throughout the section we fix a $^*-$subalgebra $\A_0 \subset L(H), \ \ 1 \in \A_0$ and an
operator $U$ such that $\d : \A_0 \to L(H)$ is a morphism (so $U$ is a partial isometry).
 Starting
from this algebra we shall extend it by means of $\d$ and $\d_*$ up to the coefficient algebra
$\A$.
\begin{I}
\label{E} \rm
 Let $X$ be a subset of $L(H)$. We denote by $\{ X \}$ the $^*-$subalgebra of $L(H)$ generated by
$X$.\\

For every $n=0,1,\ldots $ we set
\begin{equation}
\label{En}
 E_n(X) =
\{ X,\d(X), \ldots , \d^n (X)\};
\end{equation}
\begin{equation}
\label{E*n}
E_{*n}(X) = \{ X,\d_*(X), \ldots , \d_*^n (X)\}.
\end{equation}
We also set
\begin{equation}
\label{E}
E(X)= \{\ \ \bigcup_{n=0}^\infty \d^n (X)\ \  \}
\end{equation}
and
\begin{equation}
\label{E}
E_*(X)= \{\ \ \bigcup_{n=0}^\infty \d_*^n (X)\ \  \}
\end{equation}
\end{I}
\begin{prop}\ \ i) The map $\d:\A_0\to L(H)$ is a morphism iff
 $\d:E_*(\A_0)\to L(H)$ is a morphism.\\

 ii) \ \ If $\d:\A_0\to \A_0$ is an endomorphism then  $\d:E_*(\A_0)\to E_*(\A_0)$
  is an endomorphism as well.
\end{prop}
\prf \ \ i)\ \ Let us assume that $\d:\A_0\to L(H)$ is a morphism. This implies that  $U$ is a partial
isometry.
 Now, if $k\geq 1$ or $l\geq 1$ we have
\begin{equation}\d(\d_*^k(a))\d(\d_*^l(b))=U(U^{*k}aU^k)U^*U(U^{*l}bU^l)U^*=\end{equation} $$=U(U^{*k}aU^kU^{*l}bU^l)U^*=
\d(\d_*^k(a)\d_*^l(b))$$ If $k=l=0$ $\d$ is a  morphism by assumption. So, $\d:E_*(\A_0)\to L(H)$
is a morphism.

On the other hand since $\A_0$ is a subalgebra of $E_*(\A_0)$  then if $\d:E_*(\A_0)\to L(H)$ is a morphism
$\d:\A_0\to L(H)$ is a morphism as well.\\

ii) \ \ This follows from part i) and the observation that under the assumption
$\d (\A_0 )\subset \A_0$ we have for any $a\in \A_0$
$$
\d (\d_*^k (a))= UU^{*k}aU^kU^* = \d(1)\d_*^{k-1}(a)\d (1) \in E_* (\A_0).
$$
\qed
\begin{prop}
\label{A-EA*}
The following conditions
\begin{enumerate}[i)]
\item $U^*U$ is a projection and $U^*U\in\A_0'$;
\item $U^*U$ is a projection and $U^*U\in E_*(\A_0)'$;
\item $U$ and $\A_0$ have property \eqref{1};
\item $U$ and $E_*(\A_0)$ have property \eqref{1}
\end{enumerate}
are equivalent.\end{prop} \prf Equivalence\ \  i) $\Leftrightarrow$  iii)\ \ and\ \
 ii) $\Leftrightarrow$  iv)\ \ follows from Proposition \ref{equiv}.
 Clearly\ \ ii)$\Rightarrow$ i).\ \ Therefore it is enough
  to prove\ \  i)
$\Rightarrow$ ii).\ \ So let us assume that $U^*U\in\A_0'$. Then for $k\geq 0$
$$
U^*U\d_*^k(a)=U^*UU^{*k}aU^k=U^{*k}aU^k=U^{*k}aU^kU^*U=\d_*^k(a)U^*U.
 $$
  Thus $U^*U\in
E_*(\A_0)'$.
\qed
\begin{prop}
\label{E(A0)} Let $\d : \A_0 \to L(H)$ be a morphism.

 The following statements are equivalent:

(i) \ \ There exists a $^*-$algebra $\A \supset \A_0$ satisfying conditions
\eqref{1} and \eqref{one}.

(ii)
\begin{equation}
\label{2.2e} U^*U\in \bigcap_{n=0}^\infty \d^n (\A_0)'.  \end{equation} \\

%{\bf II.}\ \ If condition \eqref{2.2e} is satisfied then
%  $E(\A_0)$ is the minimal $^*-$algebra  satisfying \eqref{1} and \eqref{one} and containing
%  $\A_0$ and $\d$ is an endomorphism of $E(\A_0)$.

%  Moreover  in this case all $U^k, \ \k=1,2,\ldots$ are  partial isometries and
 % \begin{equation}
%\label{2.2einfty} U^{*k}U^k\in \bigcap_{n=0}^\infty \d^n (\A_0)', \ \ k=1,2,\ldots.
 % \end{equation}
\end{prop}
\prf \ \ (i) $\Rightarrow$ (ii). \ \ If (i) is true then as $\A \supset \A_0$ and satisfies condition
\eqref{one} we have that $\A \supset \d^n (\A_0) , \ \ n=0,1,\ldots $. This and the assumption
that $\A$ satisfies \eqref{1} implies \eqref{2.2e} in view of the equivalence
\ \ (i) $\Leftrightarrow$ (ii) \ \ from Proposition \ref{equiv}.

\ \ (ii) $\Rightarrow$     (i)\ \  Let us show that as an  algebra $\A$ satisfying
\eqref{1} and \eqref{one} one can take $E(\A_0)$.

  From \eqref{2.2e} one has for $a,b \in \A_0$ $$\d (\d^k (a)\d^l (b))= U\d^k (a)\d^l (b)U^* =
UU^*U \d^k (a)\d^l (b)U^* = U \d^k (a)U^*U\d^l (b)U^* =$$ $\d (\d^k (a))\d (\d^l (b)). $\\ Thus
$\d : E(\A_0) \to E(\A_0)$  (that is  \eqref{one} is true) and  $\d : E(\A_0) \to E(\A_0)$ is an
endomorphism.

The foregoing observation along with the note that \ \ \eqref{2.2e} $\Leftrightarrow U^*U \in
E(A_0)'$\ \ and equivalence \ \ (i)$\Leftrightarrow$ (iii)\ \ from Proposition \ref{equiv}
imply the
property \eqref{1} for $E(\A_0)$  and $U$.  \qed \\

In reality as the next proposition shows one can say much more about the algebras
satisfying  \eqref{1} and \eqref{one}.
\begin{prop}
\label{Uk}
Let $\A$ and $U$ satisfy conditions  \eqref{1} and \eqref{one} then

(i)\ \  $\d$ is an endomorphism of $\A$.

(ii)\ \
\begin{equation}
\label{k} U^k a = \d^k (a) U^k , \ \ \ \ k=1,2, ...
 \end{equation}
 for $a\in \A$.

 (iii) \ \ All the operators $U^k , \ \ k=1,2, \ldots $ are partial isometries and
\begin{equation}
\label{...}
 U^{*k}U^k \in \A', \ \ k=1,2,\ldots
\end{equation}

 (iv) The family  $U^{*k}U^k , \ \ k=1,2, \ldots $ is a commutative decreasing family of projections.

 (v) \ \  All the operators $U^{*k} , \ \ k=1,2, \ldots $ are partial isometries and
 the family  $U^kU^{*k} , \ \ k=1,2, \ldots $ is a commutative decreasing family of projections.

(vi)
\begin{equation}
\label{UkU*k-U*kUk}
[U^{*k}U^k, U^lU^{*l}] = 0,\ \ k=0,1,\ldots ,\ \ l= 0,1,\ldots
\end{equation}
where $[\alpha ,\beta ]= \alpha \beta - \beta \alpha$ is the commutator of $\alpha $ and $\beta$.

(vii) For any $1\le k \le l$ \ \
$$
 U^*U^kU^{*l} = U^{k-1}U^{*l}\ \ and \ \ UU^{*k}U^{l} = U^{*k-1}U^{l}
$$
\end{prop}
\prf
(i)\ \ Follows from \eqref{1}, \eqref{1} and the equivalence \ \ (i)$\Leftrightarrow$(iii)\ \
in Proposition \ref{equiv}.

(ii)\ \ Follows from (i) and \eqref{1}.

(iii)\ \ Follows from (ii)and the equivalence \ \ (i)$\Leftrightarrow$(ii)\ \
in Proposition \ref{equiv}.

(iv) \ \ In view of  (iii) we have that
$U^kU^{*k}U^k =U^k, \ \ k=1,2, ...$ by statement 5) of Remark
\ref{2.1}. Therefore for any $k\ge l$ we have
$$
U^{*k}U^kU^{*l}U^l = U^{*k}U^{k-l}U^lU^{*l}U^l = U^{*k}U^{k-l}U^l
= U^{*k}U^k
$$
 and in the same way one can verify that
$$
U^{*l}U^lU^{*k}U^k =  U^{*k}U^k .
$$
 Thus (iv) is true.

 (v)\ \ Since all the operators $U^k , \ \ k=1,2, \ldots $ are partial isometries  we conclude
by statement 2) of Remark
\ref{2.1}
 that
 all the
 operators  $U^{*k} , \ \ k=1,2, \ldots $ are partial isometries as well. Now by the same argument
 as used in the proof of (iv) we get (v).

 (vi) \ \ Follows from \eqref{...} and observation that $U^lU^{*l} \in \A ,\ \ l=1,2,\ldots$

 (vii)\ \ Let\ \ $1\le k\le l $\ \ then applying \eqref{UkU*k-U*kUk}
 we obtain $$
\begin{array}{*{20}l}
U^* U^k U^{*l} = (U^*U) (U^{k-1}U^{*k-1})U^{*l-k+1}=
(U^{k-1}U^{*k-1})(U^*U)U^{*l-k+1} =\\ \ \\ (U^{k-1}U^{*k-1})
(U^*UU^*) U^{*l-k} = U^{k-1}U^{*l}
\end{array}
$$
The second equality in (vii) can be proved in the same way. \qed \\

In fact Proposition \ref{Uk} and  the proof of  Proposition \ref{E(A0)} give us some
additional information which is stated in the
next proposition.
\begin{prop}
\label{Emin}
Let $\d : \A_0 \to L(H)$ be a morphism and  condition \eqref{2.2e} be satisfied then
  $E(\A_0)$ is the minimal $^*-$algebra  satisfying \eqref{1} and \eqref{one} and containing
  $\A_0$ and $\d$ is an endomorphism of $E(\A_0)$.

   Moreover  in this case all $U^k, \ \ k=1,2,\ldots$ are  partial isometries and
  \begin{equation}
\label{2.2einfty} U^{*k}U^k\in E (\A_0)'.
  \end{equation}

The families  $U^{*k}U^k , \ \ k=1,2, \ldots $ and  $U^kU^{*k} , \ \ k=1,2, \ldots $
are commutative decreasing families of projections and
\begin{equation}
\label{UkU*k-U*kUk'}
[U^{*k}U^k, U^lU^{*l}] = 0,\ \ k=0,1,\ldots ,\ \ l= 0,1,\ldots
\end{equation}
\end{prop}
\prf It has already been shown in the proof of Proposition \ref{E(A0)} that $E(\A_0)$
  satisfies \eqref{1} and \eqref{one} and that $\d$ is an endomorphism of $E(\A_0)$.
The minimality of $E(\A_0)$ is clear from the construction.

Finally as $ E(\A_0)$  satisfies   \eqref{1} and \eqref{one}  then
all the rest statements of the proposition follow from Proposition \ref{Uk}.  \qed
\begin{prop}
\label{filtrE*n}
Let $\A_0 $ possess  properties \eqref{1} and \eqref{one} then\\

(i) For any $0\le l\le k$
\begin{equation}
\label{"}
\d_*^k (\A_0)\d_*^l (\A_0)\subset \d_*^k (\A_0) \ \ {\rm and}\ \
\d_*^l (\A_0)\d_*^k (\A_0)
\subset\d_*^k (\A_0)
\end{equation}
in particular $\d_*^n (\A_0)$ is an ideal in $E_{*n}(\A_0)$ where $E_{*n}(\A_0)$ is the
algebra defined
in \ref{E}.\\

(ii) \ \ $E_{*n}(\A_0)$ is the set of operators of the form
$$
 a_0 + \d_* (a_1)+ \ldots + \d_*^N (a_n), \ \ \   a_k \in \A_0, \ \ k=0,1,\ldots , n
$$

(iii) \ \ $\d : E_{*n}(\A_0) \to E_{*n-1}(\A_0), \ \  n=1, \ldots $ is a morphism.\\

(iv)\ \  $\d_* : E_{*n}(\A_0) \to E_{*n+1}(\A_0), \ \  n=0,1, \ldots $
\end{prop}
\prf \ \ (i) \ \ By routine calculation for  $0\le l\le k$ and $a,b \in \A_0$  we have
$$
\d_*^{k}(a)\d_*^l (b) = \d_*^k (a\d^k (1)\d^{k-l}(b)) \in \d_*^k (\A_0)
$$
and
$$
\d_*^{l}(b)\d_*^k (a) = \d_*^k (\d^{k-l}(b)\d^k (1)a) \in \d_*^k (\A_0).
$$
These relations prove (i).\\
 Clearly (ii) follows from (i).\\

 (iii)\ \ Recall that $\d (\A_0)\subset \A_0 $ by the assumption
 and  for any $k\ge 1$ we have
  $$
\d (\d_*^k (a)) = UU^{*k} aU^{k}U^* = \d (1)U^{*k-1}a U^{k-1}\d (1)=
$$
$$
U^{*k-1}[\d^k(1) a \d^k (1)]U^{k-1} \in \d_*^{k-1} (\A_0) .
  $$
  This along with (ii) implies
   $\d (E_{*n}(\A_0))\subset E_{*n-1}(\A_0) , \ \ n=1,2,\ldots $. Note also that
   since $\A_0$ possesses property \eqref{1} it follows
   in view of equivalence \ \ iii) $\Leftrightarrow$ iv) \ \ from
   Proposition  \ref{A-EA*}
    that $E_{*n}(\A_0) $ being a subalgebra of  $E_*(\A_0)$ possesses this property as well.
    Therefore $\d  : E_{*n}(\A_0)\to E_{*n-1}(\A_0)$ is an endomorphism by equivalence
\ \ (i) $\Leftrightarrow$ (iii) \ \ from Proposition \ref{equiv}.\\

(iv) \ \  Evident. \qed \\

Now we are ready to present  the extensions of $\A_0$ that solve
 the problem discussed in this section : the extensions that are  coefficient algebras.
 \begin{prop}
\label{coef-ext*} Let $\A_0 $ possess  property \eqref{one}.

 The following statements are equivalent:\\

(i)\ \ There exists a coefficient algebra $\A \supset \A_0$ that is the algebra  $\A$ possessing
 properties
\eqref{1},
  \eqref{one} and  \eqref{two}.\\

  (ii) \ \ \ $\A_0$ possesses property \eqref{1}.\\
\end{prop}
\prf  \ \  (i) $\Rightarrow$ (ii). \ \ Evident as $\A_0 \subset \A$.

  (ii) $\Rightarrow$ (i)\ \  Let us show that $ E_*(\A_0)$ is a coefficient algebra.

  Since  $A_0$ possesses property \eqref{1} it follows from
 the equivalence\ \ iii) $\Leftrightarrow$ iv)\ \ in Proposition \ref{A-EA*}
   that
$ E_*(\A_0)$ possesses property \eqref{1} as well.

Now note  that in the situation under considerations all the
conditions of Proposition \ref{filtrE*n} are satisfied. Therefore
$\d : E_*(\A_0) \to E_*(\A_0)$ is an endomorphism
(in particular \eqref{one} for $E_*(\A_0)$ is satisfied) and \eqref{"} and \eqref{filtr1}
are true as well. The statement (ii) of Proposition \ref{filtrE*n} implies  the fulfillment of
\eqref{two} for $ E_*(\A_0)$.
The proof is complete. \qed \\

In fact Proposition \ref{filtrE*n} and the argument of the proof of Proposition \ref{coef-ext*}
give us some additional information which is stated in the next proposition.
\begin{prop}
\label{E*min}
 If $\A_0$ possesses properties \eqref{1} and  \eqref{one}  then $E_*(\A_0)$ is
the minimal coefficient algebra containing $A_0$ and $\d  $ is an endomorphism of $E_*(\A_0)$.

  Moreover in this situation  we have :\\
  (i) For any $0\le l\le k$
\begin{equation}
\label{"}
\d_*^k (\A_0)\d_*^l (\A_0)\subset \d_*^k (\A_0) \ \ { and}\ \
\d_*^l (\A_0)\d_*^k (\A_0)
\subset\d_*^k (\A_0).
\end{equation}
(ii) Any  element   $\beta \in E_*(\A_0)$ can be written in the form
  \begin{equation}
  \label{filtr1}
\beta = a_0 + \d_* (a_1)+ \ldots + \d_*^N (a_N)
  \end{equation}
  where $a_k \in \A_0, \ \ k=0,\ldots ,N ;\ \ N\in \N\cup\{0\}$.
\end{prop}
\prf
In the   proof of Proposition \ref{coef-ext*} it has been established that
if $\A_0$ possesses properties \eqref{1} and  \eqref{one}  then $E_*(\A_0)$ is
a coefficient algebra  and $\d  $ is an endomorphism of $E_*(\A_0)$.
The minimality of
 the extension $E_*(\A_0) \supset \A_0$
is clear from the construction.

Finally \eqref{"} and \eqref{filtr1} follow from Proposition \ref{filtrE*n}. \qed \\

As a corollary of the foregoing Propositions \ref{E(A0)}, \ref{Emin}, \ref{coef-ext*},
\ref{E*min} we obtain the following main result of this section.
\begin{thm}
\label{coef-ext} Let $\d: \A_0 \to L(H)$ be a morphism.\\

{\bf I.}\ \ The following statements are equivalent:\\

(i)\ \ There exists a coefficient algebra $\A \supset \A_0$ that is the algebra  $\A$ possessing
 properties \eqref{1},   \eqref{one} and  \eqref{two}.\\

 (ii) \ \ Condition \eqref{2.2e} is satisfied.\\

{\bf II.} \ \ If  condition \eqref{2.2e} is satisfied then $E_*(E(\A_0))$ is the minimal
coefficient algebra containing $\A_0$ and $\d$ is an endomorphism of $E_*(E(\A_0))$.

 Moreover we have\\
  (i)\ \ For any $0\le l\le k$
$$
\d_*^k (E(\A_0))\d_*^l (E(\A_0))\subset \d_*^k (E(\A_0)) \ \ { and}\ \
\d_*^l (E(\A_0))\d_*^k (E(\A_0))
\subset\d_*^k (E(\A_0)).
$$
(ii)\ \ Any  element   $\beta \in E_*(E(\A_0))$ can be written in the form
  $$
\beta = \alpha_0 + \d_* (\alpha_1)+ \ldots + \d_*^N (\alpha_N)
  $$
  where $\alpha_k \in E(\A_0), \ \ k=0,\ldots ,N ;\ \ N\in \N\cup\{0\}$.\\
  (iii)\ \  All\ \ $U^k, \ \ k=1,2,\ldots$ are  partial isometries and
  \begin{equation}
\label{2.2einfty} U^{*k}U^k\in E_*(E (\A_0))'.
  \end{equation}

The families  $U^{*k}U^k , \ \ k=1,2, \ldots $ and  $U^kU^{*k} , \ \ k=1,2, \ldots $
are commutative decreasing families of projections and
\begin{equation}
\label{UkU*k-U*kUk'}
[U^{*k}U^k, U^lU^{*l}] = 0,\ \ k=0,1,\ldots ,\ \ l= 0,1,\ldots
\end{equation}
\end{thm}

Summarizing the extension procedure $\A_0 \hookrightarrow  E_*(E(\A_0)) $ described above we
conclude that if we have a $^*-$algebra $\A_0 \subset L(H), \ \ 1\in A_0$  such that
$\d : \A_0 \to L(H)$ is a morphism then
the $C^*-$algebra $\B = \B(\A_0, U)$ generated by $\A_0$ and $U$ could be equivalently given as
$\B =\B (E_*(E(\A_0)), U)$ where $E_*(E(\A_0))$ is the minimal coefficient algebra containing $\A_0$.

\section{Commutative coefficient algebra}

Considering the extension procedures $\A_0 \hookrightarrow E(\A_0)$, \ \
$\A_0 \hookrightarrow E_*(\A_0)$\ \ and \ \ $\A_0 \hookrightarrow   E_*(E(\A_0))$ presented in the
previous section one naturally arrives at the question: what properties of $\A_0$ are preserved
under these procedures? It is clear that to preserve a certain property, for example,
 commutativity  one needs to take into account some additional interrelations between
 $\A_0$ and $U$. Precisely this subject is the theme of the present section.

 Throughout this section we assume that $\A_0$ is a
 commutative $^*-$subalgebra of $L(H)$ and
 $1\in \A_0$.
\begin{prop}
\label{commE*}
Let $\A_0$ be a commutative $^*-$subalgebra of $L(H)$  satisfying conditions
\eqref{1} and \eqref{one} then $E_* (\A_0)$ is a commutative  algebra and both the mappings
$\d : E_* (\A_0) \to E_* (\A_0)$ and $\d_* : E_* (\A_0) \to E_* (\A_0)$ are endomorphisms.
\end{prop}
\prf In the situation under consideration all the assumptions of Proposition \ref{filtrE*n}
are satisfied.
The argument in the proof of (i) of Proposition \ref{filtrE*n} along with the commutativity
of $\A_0$  and the condition  $\d (\A_0) \subset \A_0$ implies the equality
$$
\d_*^k(a)\d_*^l (b) =  \d_*^l (b)\d_*^k(a), \ \ k,l \ge 0,\ \ a,b,\in \A_0 .
$$
Thus  $E_* (\A_0)$ is commutative.

By Proposition \ref{E*min}\ \ $\d : E_* (\A_0) \to E_* (\A_0)$ is an endomorphism
and $\d_* ( E_* (\A_0)) \subset  E_* (\A_0)$ so to finish the
proof it is enough to verify that $\d_* : E_* (\A_0) \to E_* (\A_0)$ is an endomorphism.

Since $E_* (\A_0)$ is commutative and $\d (\A_0)\subset \A_0$ we have that
$$
UU^* = \d(1) \in \A_0 \subset [E_* (\A_0)]'.
$$
In view of this we can apply the equivalence  \ \ (ii)$\Leftrightarrow$(iii)\ \ in Proposition
\ref{equiv} to $U^*$, $\d_*$ and
$E_* (\A_0)$ and conclude that $\d_* : E_* (\A_0) \to E_* (\A_0)$ is an endomorphism. \qed \\
 \begin{thm}
 \label{commE-E*A}
Let $\A_0$ be a commutative $^*-$subalgebra of $L(H)$ and  $\d : \A_0 \to L(H)$ be a morphism.\\

The following three statements are equivalent:\\

(i)\ \ There exists a commutative $^*-$algebra $\overline{\A} \supset \A_0$
satisfying conditions \eqref{1} and \eqref{one}.\\

(ii)\ \ There exists a commutative coefficient algebra $\A \supset \A_0$.\\

(iii)\ \ The following two conditions
\begin{equation}
\label{commEA}
\A_0 \subset \bigcap_{n=0}^\infty \d^n (\A_0)'
\end{equation}
and
\begin{equation}
\label{2.2e'} U^*U\in \bigcap_{n=0}^\infty \d^n (\A_0)'.
\end{equation}
are satisfied.
\end{thm}

 \prf \ \ (i)$\Rightarrow$ (iii)\ \ Here \eqref{2.2e'} is true by Proposition \ref{E(A0)}.
 In addition since $\A_0 \subset \overline{\A} $ and $\overline{\A}$ is commutative and
 satisfies  \eqref{one} it follows that \eqref{commEA} is true.\\

 (ii)$\Rightarrow$ (iii)\ \ Here \eqref{2.2e'} follows from Theorem \ref{coef-ext}.
  In addition since $\A_0 \subset {\A} $ and ${\A}$ is commutative and
 satisfies  \eqref{one} it follows that \eqref{commEA} is true.\\

 (iii)$\Rightarrow$ (i)\ \  Let us take here  $\overline{\A} = E(\A_0)$.
      By Proposition \ref{Emin} $E(\A_0)$ satisfies \eqref{1} and \eqref{one} so it is enough to
      verify the commutativity of  $E(\A_0)$.

      In view of \eqref{commEA} we have
      \begin{equation}
 \label{000}
[a,\d^n (b)] = 0, \ \ a,b\in \A_0, \ \ n=0,1,\ldots
 \end{equation}
  In addition
   \eqref{2.2e'} and Proposition \ref{Emin} imply that
 $\d : E(\A_0) \to E(\A_0)$
 is a morphism. This along with \eqref{000} gives
 $$
[\d^{k+l}(a), \d^l (b)] = \d^l ([\d^k (a),b]) = 0.
 $$
 Thus $E(\A_0)$ is commutative.\\

(iii)$\Rightarrow$ (ii)\ \  Let us take here $\A = E_* (E(A_0))$. By \eqref{2.2e'}
and Theorem \ref{coef-ext} $ E_* (E(A_0))$ is a coefficient algebra so it is enough
 to prove that  $ E_* (E(A_0))$ is commutative. But we have already proved that
$E(\A_0)$ is commutative so by Proposition \ref{commE*}  $ E_* (E(A_0))$ is
commutative as well. \qed \\

In fact we can strengthen the foregoing result.
\begin{prop}
\label{E*E-comm-min}
Let $\A_0$ be a commutative $^*-$subalgebra of $L(H)$ and  $\d : \A_0 \to L(H)$ be a morphism.
  If \eqref{commEA} and \eqref{2.2e'} are satisfied then $E(\A_0)$ is the
 minimal commutative $^*-$algebra containing $\A_0$ and satisfying \eqref{1} and \eqref{one}, and
 $E_* (E(\A_0))$ is the minimal commutative coefficient algebra containing $\A_0$.
 Moreover $\d : E_* (E(\A_0)) \to E_* (E(\A_0))$ and $\d_* : E_* (E(\A_0)) \to E_* (E(\A_0))$
 are endomorphisms and
 \begin{equation}
 \label{EE*=E*E}
  E_* (E(\A_0))= E (E_*(\A_0))
  \end{equation}
\end{prop}
\prf
By the argument of the proof of Theorem \ref{commE-E*A} $E(\A_0)$
     satisfies \eqref{1} and \eqref{one}  and $E_* (E(\A_0))$ is a coefficient algebra.
The minimality of $E(\A_0)$ and  $E_* (E(\A_0))$  follows from the construction.

Since $E(\A_0)$ satisfies \eqref{1} and \eqref{one} Proposition \ref{commE*} tells us that
$\d : E_* (E(\A_0)) \to E_* (E(\A_0))$ and $\d_* : E_* (E(\A_0)) \to E_* (E(\A_0))$
 are endomorphisms. So to finish the proof we have to verify \eqref{EE*=E*E}.

 Observe that
\begin{equation}
 \label{???}
\d_* : \A_0 \to E_* (E(\A_0))\ \ {\rm is \ \ a \ \ morphism},
 \end{equation}
 \begin{equation}
 \label{?}
UU^* \in \A_0 \subset E_* (E(\A_0)) \subset \bigcap_{n=0}^\infty \d_*^n (\A_0)',
 \end{equation}
 and
 \begin{equation}
 \label{??}
 \A_0 \subset E_* (E(\A_0)) \subset \bigcap_{n=0}^\infty \d_*^n (\A_0)'.
 \end{equation}
 Note that\eqref{?} means that we have substituted $U$ by $U^*$ and $\d$ by $\d_*$ in
  \eqref{2.2e'}.\ \ On the other hand
\eqref{??}   means that we have substituted  $\d$ by $\d_*$ in \eqref{commEA}.
 Thus by the already
proved part of the proposition both the mappings
$\d : E (E_*(\A_0)) \to E (E_*(\A_0))$ and $\d_* : E (E_*(\A_0)) \to E (E_*(\A_0))$
are endomorphisms. This implies
\begin{equation}
\label{E*EE*}
E_*(E (E_*(\A_0))) \subset E (E_*(\A_0)).
\end{equation}
But since
$$
E_*(E (E_*(\A_0))) \supset  E_* (E(\A_0))
$$
\eqref{E*EE*} means that
\begin{equation}
\label{E*EE*'}
E_*(E (\A_0)) \subset E (E_*(\A_0)).
\end{equation}
By the same argument we have
\begin{equation}
\label{E*EE*''}
E (E_*(\A_0)) \subset E_*(E (\A_0)).
\end{equation}
Thus
\eqref{EE*=E*E} is true.

\section{Coefficient algebras generated by polar decomposition}

In this section we shall present  interesting examples of  coefficient algebras. As
we shall see these examples are important for many reasons in particular  for application
in the theory of quantum physical systems.

{\bf Example 1.} \ \
Let $a\in L(H)$ be a certain operator and
\begin{equation}
 \label{polar}a=U\abs a
 \end{equation}
be the standard polar decomposition
of $a$. Here $\abs a = \sqrt {a^* a}$
and   $U$   is
a partial isometry
defined by
\begin{equation}
\label{eI.3}
U(|a|\xi ) = a\xi , \ \ \ \xi \in H .
\end{equation}
As the algebra $\A_0$ we take
   the commutative $C^*-$algebra
\begin{equation}
 \A_0=\overline{ \{1,\abs a\}}.
 \end{equation}
Everywhere further in this example   we shall assume the satisfaction of the condition
\begin{equation}
 \label{aa}
 aa^*\in\A_0,
\end{equation}
{\bf Remark.} $C^*-$algebras defined by condition \eqref{aa} are natural generalizations
of the algebras satisfying the relation
\begin{equation} aa^*=\gamma(a^*a),\end{equation}
where $\gamma$ is a continuous positive-valued function on the spectrum of the  operator $a^*a$.

Note that since in the situation under consideration  $U^*U$ is the orthogonal
projection onto ${\rm Im}\, \abs a = ({\rm Ker} \abs a )^\bot$  it follows that
$U^*U$ is the spectral projection of $\abs a$ corresponding to the set
\ \ $\sigma (|a|) \setminus \{ 0 \}$ \ \ (where we denote by $\sigma
(\alpha )$ the spectrum of an operator $\alpha$). Thus by the
spectral theorem
 (see, for example, \cite{Naimark}, \S 17)
\begin{equation}
\label{bicommutant}
U^*U \in \overline{\{ 1, |a| \}}^{\prime \prime} = \A_0''.
\end{equation}
where $\A_0''$ denotes the bicommutant of $\A_0$ (that is the Von Neumann algebra
 generated by $\A_0$). In particular we have
\begin{equation}
 U^*U\in\A_0'.
\end{equation}
So according to Proposition \ref{equiv}
\begin{equation}
 \d:\A_0\to \L(\H)
 \end{equation}
is a morphism and $U$ and $A_0$ satisfy  condition \eqref{1}.

The next theorem  describes  a number of  properties of $\A_0$ and $U$ that follow from
condition \eqref{aa}.
\begin{thm}
\label{polarU}
Let $U$ be the partial isometry defined by polar decomposition \eqref{polar} of an operator $a$
satisfying \eqref{aa}. Then:
\begin{enumerate}[(i)]
\item
\begin{equation}
 \d^k(\{\abs a\})\subset \overline{ \{1,\abs a,UU^*,\ldots,U^{k-1}U^{*k-1} \}},\ \ k=1,2, ...
 \end{equation}
\item
$\d$ is an endomorphism of the $C^*-algebra$
 $$
 \overline{ \{1,\abs a,UU^*,\ldots,U^{k}U^{*k},\ldots ;\ \ k=1,2, \ldots  \}}
 $$
\item
\begin{equation}
 U^kU^{*k}=\d^k(1)\in\A_0'', \ \ k=0,1,\ldots
 \end{equation}
\item
 \label{t4}
 For $1\leq k\leq l$
  \begin{equation} U^*U^kU^{*l}=U^{k-1}U^{*l}\ \ and \ \  UU^{*k}U^{l}=U^{*k-1}U^{l}
  \end{equation}
\item \begin{equation} [U^{*l}U^l,U^kU^{*k}]=0
\end{equation}
for $k,l\in\N$.
\end{enumerate}\end{thm}
\prf
(i)\ \ Observe that property \eqref{aa} means that
\begin{equation}
\label{i}
\d (\abs{a}^2)= U\abs{a}^2U^* \in \overline{\{ 1, \abs a  \}}
\end{equation}
Since $\d : \A_0 \to L(H)$ is a morphism it follows that
$$
  [\d (\abs a) ]^2 =\d (\abs{a}^2)\in \overline{\{ 1, \abs a  \}}.
$$
Therefore \eqref{i} implies
\begin{equation}
\label{i'}
\d (\abs{a})\in \overline{\{ 1, \abs a  \}}.
\end{equation}
By iteration of \eqref{i'} we have (i).

(ii)\ \ Follows from (i) and the fact that  $\d : \A_0 \to L(H)$ is a morphism.

(iii)\ \  Since as it has been already observed  $U^*U$ is the spectral projection
corresponding to the interval $(0, |a|]$  it
follows that there exists a sequence $\alpha_n$ of elements of $\{
|a|  \}$ such that
\begin{equation}
\label{strong}
\alpha_n\ \  \mathop{\longrightarrow}\limits^{strongly}\ \ U^*U
\end{equation}
Therefore we have $$ U\alpha_n U^* \ \
\mathop{\longrightarrow}\limits^{strongly}\ \ U (U^*U)U^* = UU^*
$$ Since (by (i)) $U\alpha_n U^* \subset \overline {\{ 1,|a| \}}$ it follows
that $UU^* \in \overline{ \{ 1,|a| \}}^{\prime\prime}$.\\ The further
proof goes by induction.\\ Suppose that\ \
$U^kU^{k*}\in  \overline {\{ 1,|a| \}}^{\prime\prime} , \ \ k=\overline{1, n-1}$. Taking the
sequence $\alpha_n$ mentioned above we have
\begin{equation}
\label{e2.2.1}
U^n\alpha_n U^{n*} \ \  \mathop{\longrightarrow}\limits^{strongly}\ \ U^{n} (U^*U)U^{n*} =
U^{n-1} (UU^*U) U^{n*} =
 U^{n}U^{n*}
\end{equation}
But due to (i) and the assumption of the induction we have $$
U^n\alpha_n U^{n*} \subset \{ 1, |a|, UU^*, ... , U^{n-1}U^{n-1*}
\}\subset \overline{ \{ 1,|a| \}}^{\prime\prime} $$ and therefore
(\ref{e2.2.1}) implies \ \
 $U^{n}U^{n*} \in \overline{ \{ 1,|a| \}}^{\prime\prime}$. So (iii) is proved.

 (iv), (v)\ \ By (iii) we have that
 $$
\overline{\A} := \overline{ \{1,\abs a,UU^*,\ldots,U^{k}U^{*k},\ \ k=1,2, \ldots  \}}
\subset \A_0^{''}
 $$
and in particular $\overline{\A}$ is a commutative algebra.

In view of (ii) $\d :\overline{\A} \to \overline{\A}$ is an endomorphism
and it has already been observed that
$$
U^*U \in \A_0'' \subset \overline{\A}' .
$$
Bearing in mind the equivalence \ \ (i)$\Leftrightarrow$(ii)\ \ in Proposition \ref{equiv}
we conclude that $\overline{\A} $ and $U$ satisfy all the assumptions of
Proposition \ref{Uk}. This implies (iv) and (v). \qed

\begin{rem}
\rm
Most of the properties listed in Theorem \ref{polarU} are known
(we have presented them for the sake of completeness). In particular one can find some generalizations
of property (v) and also generalizations of particular cases of (i) in Propositions 28 and 29 in
\cite{OstSam}, Section 2.1 that also contains a lot of important information related to the subject
considered.
\end{rem}

As the corollary of the properties listed  and the results obtained in the previous section
 we can get the following result.
\begin{prop}
\label{polar1}
Let $U$ be the partial isometry defined by polar decomposition \eqref{polar} of an operator $a$
satisfying \eqref{aa} and $\A_0 = \overline{ \{1,\abs a \}}$ then

(i)\ \
\begin{equation}
\label{EApolar}
E(\A_0) = \{\A_0 , UU^*,\ldots,U^{k}U^{*k},\ldots   \}
\end{equation}
and $E(\A_0)$ is a commutative $^*-$algebra satisfying conditions \eqref{1} and \eqref{one}
and $\d : E(\A_0) \to E(\A_0)$ is an endomorphism.

(ii)\ \ The algebra $E_*(E(\A_0))$ is a commutative coefficient subalgebra and
both the mappings
 $\d : E_* (E(\A_0)) \to E_* (E(\A_0))$ and $\d_* : E_* (E(\A_0)) \to E_* (E(\A_0))$
 are endomorphisms.
\end{prop}
\prf \ \
  The statements (i), (iii) of Theorem \ref{polarU} and \eqref{bicommutant}
imply the satisfaction of  properties \eqref{commEA}        and     \eqref{2.2e'}.
Therefore Proposition \ref{E*E-comm-min} tells us that
 $E(\A_0)$ is a commutative $^*-$algebra  satisfying \eqref{1}
 and \eqref{one} and $\d : E(\A_0) \to E(\A_0)$ is an endomorphism; and
 $E_* (E(\A_0))$ is a commutative coefficient algebra and
  the mappings
   $\d : E_* (E(\A_0)) \to E_* (E(\A_0))$ and $\d_* : E_* (E(\A_0)) \to E_* (E(\A_0))$
 are endomorphisms.

 The equality \eqref{EApolar} follows from statement (i) of Theorem \ref{polarU}. \qed

 {\bf Example 2.} \ \ Let us consider a partial isometry $U\in L(H)$ and a positive operator
 $0\le Q \in L(H)$ satisfying the conditions
 \begin{equation}
 \label{Q1}
U^*U = \d_* (1)\in \{ Q \}'' \subset \A_0''
 \end{equation}
 and
 \begin{equation}
 \label{Q2}
\d (Q) = UQU^* = qQ
 \end{equation}
 where
 $\A_0 = \{1,Q  \}$ and $0<q<1$.

 By Proposition \ref{equiv} $\d : \A_0 \to L(H)$ is a morphism and in addition from \eqref{Q2}
  it follows that
  \begin{equation}
  \label{Q3}
\d (f (Q)) = f (qQ)
  \end{equation}
  for any $f\in C(\sigma )$  such that $f(0) = 0$ (here $ C(\sigma )$  is the algebra of continuous functions on the
  spectrum $\sigma $ of operator $Q$). From \eqref{Q3} we conclude that $\d$ maps the algebra
  $\{  Q \}$ on itself and therefore the condition
  \begin{equation}
  \label{QQQQ}
U^*U \in \bigcap_{n=0}^\infty \d^n (\A_0 )'
  \end{equation}
  is fulfilled iff
  \begin{equation}
  \label{QQ}
[U^*U, \d^n (1)] =0, \ \ n=1,2, \ldots
  \end{equation}
Following the argument used in the proof of statement (iii) of Theorem \ref{polarU}
and using
\eqref{Q1} instead of \eqref{bicommutant} one can prove that
\begin{equation}
\label{bic''}
\d^n(1)\in\A_0'', \ \ \ n=1,2,\ldots
\end{equation}
This along with \eqref{Q1} proves \eqref{QQ} and \eqref{QQQQ}.

Similarly the condition
$$
\A_0 \subset \bigcap_{n=0}^\infty \d^n (\A_0)'
$$
can be reduced to the condition
$$
[a, \d^n (1)]=0, \ \ n=1,2, \ldots ; a\in \A_0
$$
which is again satisfied  in view of   \eqref{bic''}.

Therefore Proposition \ref{E*E-comm-min} tells us that $E_* (E(\A_0))$ is
a coefficient algebra. Hence one can use the isomorphism Theorem \ref{iso} to investigate the
structure of the algebra $\B = \B (\A_0, U)$.

Let us consider the operator $a \in \B \subset L(H)$ given by the formula
\begin{equation}
\label{Q4}
a:= U \rho (Q)
\end{equation}
where $0\le \rho \in C(\sigma )$ and $\rho (0)=0, \ \ \rho (t)\neq 0, \ \ t\neq 0$
 and $U^*U \rho (Q)= \rho (Q)$.

 In view of \eqref{Q1} we have that \eqref{Q4}
is the polar decomposition of $a$.

\begin{prop}
\label{relations}
Operators $a,a^*$ and $Q$ satisfy the following relations:
\begin{equation}
\label{Q5}
a^*a =\rho ^2 (Q),
\end{equation}
\begin{equation}
\label{Q6}
aa^* = \rho ^2 (qQ),
\end{equation}
\begin{equation}
\label{Q7}
aQ = qQa,
\end{equation}
\begin{equation}
\label{Q8}
Qa = qa^*Q,
\end{equation}
\end{prop}
\prf
As by definition \eqref{Q4} is the polar decomposition then \eqref{Q5} is true.

 Equality \eqref{Q6} follows from \eqref{Q3}.

To verify \eqref{Q7} note that property \eqref{Q2} implies
\begin{equation}
\label{Q9}
Ub =\d (b)U, \ \ b\in \A_0.
\end{equation}
Substituting here $b$ by $Q$ and using \eqref{Q2} one obtains
\begin{equation}
\label{Q10}
UQ =qQU
\end{equation}
that gives
\begin{equation}
\label{Q11}
U\rho(Q)= \rho(qQ)U.
\end{equation}
From \eqref{Q11} one has
$$
aQ = U\rho (Q)Q = \rho (qQ)qQU = qQU\rho (Q) = qQa .
$$
Thus \eqref{Q7} is true.

Finally \eqref{Q8} follows from the previous formula by conjugation. \qed \\

The $C^*-$algebra $\B_\rho := \B (1,a,Q)$ generated by the operators $1,a$ and $Q$
satisfying \eqref{Q5}-\eqref{Q8} is a $C^*-$subalgebra of $\B (\A_0, U)$.

In a special case when
\begin{equation}
\label{Q12}
\rho (Q) = \frac{1-Q}{1-q}
\end{equation}
the algebra $\B_\rho$ will be the $q-$deformation Heisenberg algebra (see \cite{MaxOdz}).
And in  the case when
\begin{equation}
\label{Q13}
\rho ^2 (Q) = - \frac{1}{(1-q)(1-q^2)}\left( \frac{q}{Q}+Q \right) \frac{1}{(1-q)^2}
\end{equation}
$\B_\rho $ is the
 $q-$deformation $U_q ({\rm Sl} (2))$ enveloping algebra of ${\rm Sl} (2)$ (see \cite{Odz1}).

 The general case of $\B_\rho$ was studied  in \cite{Odz1} where an application of this object to the
 theory of basic hypergeometric series and integrable physical systems was investigated.

 The algebras $\B_\rho$ also appear in a natural way in quantum optics, namely, by the
 quantum reduction of the  multi-mode quantum optical system. Reducing
 the degree of freedom of
 the quantum system considered one passes from the multi freedom degree   Heisenberg
 algebra to the one given by relations \eqref{Q5}-\eqref{Q8}. The
  operators $a$ and $a^*$ have the interpretation of the cluster annihilation and creation
  operators. For an exhaustive description of this important physical
   application see \cite{OHT2}.

  \end{document}